\documentstyle[fleqn,12pt]{article}
\begin{document}
\pagenumbering{arabic}
\setcounter{page}{1}
\pagestyle{plain}
\baselineskip=20pt

\thispagestyle{empty}
\rightline{YTUMB 00-01, May 2001} 
\vspace{1.4cm}

\begin{center}
{\Large\bf Differential Geometry of the $q$-Quaternions}
\end{center}

\vspace{1cm}
\begin{center} Salih \c Celik\footnote{E-mail: sacelik@yildiz.edu.tr}\\
Yildiz Technical University, Department of Mathematics, \\
34210 Davutpasa, Istanbul, TURKEY. \end{center}

\vspace{2cm}
{\bf Abstract}

Differential calculus on the quantum quaternionic  group GL$(1, H_q)$ is 
introduced. 

\vfill\eject
{\bf 1. Introduction}

Differential geometry in the theory of (quantum) Lie groups plays 
an important role in the mathematical modelling of physics theories. 
In the classical differential geometry one has a choice between two dual 
and equivalent descriptions: one can either work with points on a manifold 
$M$ or with the algebra $C(M)$ of smooth functions on 
$M$. The idea that the algebra $C(M)$ need not be commutative 
gives rise to the noncommutative geometry. Such a space is called a quantum 
space. This is analogy with the quantization of the commutative algebra of functions 
on phase space that yields the noncommutative operator algebra of quantum 
mechanics. 

A class of noncommutative Hopf algebras have been found in the discussions 
of integrable systems. These Hopf algebras are $q$-deformed function 
algebras of classical groups and this structure is called quantum group [1]. 
The quantum group can also be regarded as a generalization of the notion of 
a group [2]. Noncommutative geometry [3] is one of the most attractive 
mathematical concepts in physics and has started to play an important role 
in different fields of mathematical physics for the last few years. 
The basic structure giving a direction to the noncommutative geometry 
is a differential calculus [4] on an associative algebra. 

Quantum quaternionic algebra and its Hopf algebra structure is important 
in physics. The importance of differential geometry in the quantum 
quaternionic algebra should not be underestimated. Quantum quaternion 
is an example of a quantum space, and to investigate its differential 
geometry may be interesting. This is considered in the present work. 

\vfill\eject
{\bf 2. Review of Hopf algebra $H_q$ }

Elementary properties of quantum quaternionic group GL$(1,H_q) = H_q$ are 
described in Refs. 5 and 6. We state briefly the properties we are going to 
need in this work. 

{\bf 2.1 The algebra of functions on $H_q$ }

The quantum quaternionic algebra $H_q$ is defined as a pair $({\cal A, M})$ 
equipped with $*$ structure, where ${\cal A}$ is an algebra and ${\cal M}$ 
is an ${\cal A}$-module, and they have the following properties: 

{\bf (1)} ${\cal A}$ is an unital associative algebra generated by generators 
$a_k$ $(k = 0, 1, 2, 3)$ with the commutation relations 
$$a_0 a_1 = a_1 a_0 - {\bf i} {{q - q^{-1}}\over 2} (a^2_2 + a^2_3), $$
$$a_0 a_2 = {{q + q^{-1}}\over 2} a_2 a_0 + 
            {\bf i} {{q - q^{-1}}\over 2} a_2 a_1, $$
$$a_0 a_3 = {{q + q^{-1}}\over 2} a_3 a_0 + 
            {\bf i} {{q - q^{-1}}\over 2} a_3 a_1, \eqno(1)$$
$$a_1 a_2 = {{q + q^{-1}}\over 2} a_2 a_1 - 
            {\bf i} {{q - q^{-1}}\over 2} a_2 a_0, $$
$$a_1 a_3 = {{q + q^{-1}}\over 2} a_3 a_1 - 
            {\bf i} {{q - q^{-1}}\over 2} a_3 a_0, $$
$$a_2 a_3 = a_3 a_2, $$
where ${\bf i}^2 = - 1$ and $q$ is a nonzero real number. 

{\bf (2)} The $*$ antiinvolution in ${\cal A}$ is defined by 
$$ a_{j}^{\star} = a_j, \qquad j = 0, 1 $$
$$ a_{2}^{\star} = {1\over 2}\left[(q + q^{-1}) a_2 - {\bf i} 
                    (q - q^{-1}) a_3\right), 
\eqno(2) $$
$$a_{3}^{\star} = {{{\bf i}}\over 2}\left[(q - q^{-1}) a_2 - {\bf i} 
                   (q + q^{-1}) a_3\right]. $$
Note that
$$ (a_{k}^{\star})^{\star} = a_k, \qquad k = 0, 1, 2, 3. $$

{\bf (3)} ${\cal M}$ is an ${\cal A}$-module generated by the quaternionic 
units $e_k$ with the relations 
$$ e_k e_l = - \delta_{kl} e_0 + \epsilon_{klm} e_m \eqno(3) $$
and 
$$a_k e_l = e_l a_k, \eqno(4)$$
where $\delta_{kl}$ denotes the Kronecker delta and 
$$ \epsilon_{klm} = {1\over 2}(k - l)(l - m)(m - k). $$

The quaternionic conjugation ~($^\star$-antiinvolution) in ${\cal M}$ 
is defined by
$$ e _{k}^{\star} = 2 \delta_{0,k}e_0 - e_k, \qquad k = 0,1,2,3. \eqno(5)$$

{\bf (4)} Assume that any quaternion $h$ has a representation 
$$ h = a_0e_0 + a_1e_1 + a_2e_2 + a_3e_3 \eqno(6) $$
in terms of the generators of ${\cal A}$. $h$ will be called the $q$-quaternion, 
and in this case we shall say that $q$-quaternion $h$ belongs to $H_q$. 

The conjugation of a $q-$quaternion $h$ is introduced as 
$$ h^{\star} = e_0 a_{0}^{\star} - e_1 a_{1}^{\star} - e_2 a_{2}^{\star} - 
               e_3 a_{3}^{\star}. \eqno(7) $$
Hence, we can introduce the $q$-norm of $h$ as 
$$ {\cal N}_q(h) = hh^{\star} = 
      a_{0}^{2} + a_{1}^{2} + {1\over 2} (q + q^{-1}) (a_{2}^{2} + a_{3}^{2}).
\eqno(8) $$
Note that ${\cal N}_q(h)$ belongs to the center of $H_q$. 

{\bf 2. 2 The Hopf Algebra Structure of $H_q$ }

The action of comultiplication $\Delta$ on the generators $a_k$ of ${\cal A}$ 
can be introduced as 
$$ \Delta(a_0) = a_0 \otimes a_0 - (a_1 \otimes a_1 + a_2 \otimes a_2 + 
a_3 \otimes a_3), $$
$$ \Delta(a_1) = a_0 \otimes a_1 + a_1 \otimes a_0 + a_2 \otimes a_3 - 
a_3 \otimes a_2 , \eqno(9) $$
$$ \Delta(a_2) = a_0 \otimes a_2 + a_2 \otimes a_0 + a_3 \otimes a_1 - 
a_1 \otimes a_3 , $$
$$ \Delta(a_3) = a_0 \otimes a_3 + a_3 \otimes a_0 + a_1 \otimes a_2 - 
a_2 \otimes a_1. $$
Note that 
$$ \Delta(e_0) = e_0 \otimes e_0. \eqno(10) $$
It also easy to show that 
$$ \Delta({\cal N}_q(h)) = {\cal N}_q(h) \otimes {\cal N}_q(h). \eqno(11) $$

The action of counit $\varepsilon$ on the generators $a_k$ of ${\cal A}$ 
is given by
$$ \varepsilon(a_k) = \delta_{0,k} e_0, \qquad k = 0, 1, 2, 3 \eqno(12) $$
$$ \varepsilon(e_0) = e_0. $$

The action of antipode ${\cal S}$ on the generators $a_k$ of ${\cal A}$ 
is introduced as 
$$ {\cal S}(a_k) = {\cal N}_{q}^{-1}(h)(2 \delta_{0,k}a_0 - a_{k}^{\star}) 
\eqno(13) $$
for $k = 0, 1, 2, 3$. 

{\bf 3. Differential calculus on $H_q$} 

In this section, we shall build up the differential calculus on the 
quantum quaternionic algebra $H_q$. The differential calculus on 
$H_q$ involves functions on the algebra, differentials and differential forms. 

{\bf 3.1 Classical case}

Let's begin with differential calculus on the classical quaternionic group 
GL$(1,H)$. In a classical Lie group $G$, one-order differential calculus is a 
linear map 
$${\sf d} : C^\infty(G) \longrightarrow \Gamma,$$
where $C^\infty(G)$ is a ${\cal C}$-algebra consisting of all smooth functions 
on $G$ and $\Gamma$ is a $C^\infty(G)$-bimodule consisting of all differential 
forms. This linear map satisfies 

\noindent
{\bf (i)} the nilpotency 
$${\sf d}^2 = 0, \eqno(14)$$
{\bf (ii)} for all $f, g \in C^\infty(G)$, 
$${\sf d} (f g) = ({\sf d} f) ~g + (-1)^{\hat{f}} f ~({\sf d} g) \eqno(15)$$
where $f$ and $g$ are functions of the generators and $\hat{f}$ is the 
corresponding grading of $f$. 

According to the ideas of noncommutative geometry [3], the differential 
calculus can be defined on a more general noncommutative algebra. 

Let ${\cal A}_1$ be a Hopf algebra with unit generated by the generators 
$a_k$ $(k = 0,1,2,3)$. We denote differentials of $a_k$ by ${\sf d} a_k$. 
Then one can construct the one-form $\Omega$, where 
$$\Omega = {\sf d} h ~h^\star. \eqno(16)$$
Explicitly, 
$$w_0 = {\sf d} a_0 a_0 + {\sf d} a_1 a_1 + {\sf d} a_2 a_2 + 
  {\sf d} a_3 a_3, \eqno(17)$$
etc. It can be easily checked that these one-forms construct a four 
dimensional Grassmann algebra. The anti-commutation relations of the one-forms 
allow us to construct the algebra of the generators. To obtain the Lie algebra 
of the algebra generators we first write the one-forms as
$${\sf d} a_0 = w_0 a_0 - w_1 a_1 - w_2 a_2 - w_3 a_3, \eqno(18)$$
etc. The differential {\sf d} can then the expressed in the form 
$${\sf d} = 2 (w_0 \nabla_0 - w_1 \nabla_1 + w_2 \nabla_2 - w_3 \nabla_3),
\eqno(19)$$ 
where $\nabla_k$ $(k=0,1,2,3)$ are the Lie algebra generators. We wish to 
obtain the commutation relations of these generators.  Let $f$ be an arbitrary 
function of the generators of ${\cal A}_1$. Then, using the nilpotency of 
the exterior differential {\sf d} we obtain 
$$(-1)^k {\sf d} w_k ~\nabla_k f = (-1)^{k+j} w_k (w_j \nabla_j) \nabla_k f, 
\qquad k, j = 0,1,2,3 \eqno(20)$$
where summation over repeated indices is understood. Using one-forms one 
easily obtain the two-forms 
$${\sf d} w_0 = 0, \qquad {\sf d} w_1 = 2 w_2 w_3, $$
$${\sf d} w_2 = 2 w_3 w_1, \qquad {\sf d} w_3 = - 2 w_2 w_1, \eqno(21)$$
since 
$${\sf d} \Omega = - {\sf d} h ~{\sf d} h^\star = \Omega^2. \eqno(22)$$
We now find the following commutation relations for the Lie algebra 
$$[\nabla_1,\nabla_0] = [\nabla_2,\nabla_0] = [\nabla_3,\nabla_0] = 0,$$
$$[\nabla_1,\nabla_2] = - 2 \nabla_3, \qquad [\nabla_2,\nabla_3] = - 2 \nabla_1, 
  \qquad [\nabla_3,\nabla_1] = - 2 \nabla_2. \eqno(23)$$

{\bf 3.2 Quantum case} 

A differential algebra on $H_q$ is an associative algebra 
$\Gamma$ equipped with an operator {\sf d}. Also the algebra $\Gamma$ 
has to be generated by ${\cal A} \cup {\sf d} {\cal A}$. 

Firstly, to obtain the relations between the generators of ${\cal A}$ 
and their differentials. We shall use the method of Ref. 7. Using the 
consistency of a differential calculus, as the final result one has 
the following commutation relations 
$$a_0 ~{\sf d} x_+ = {{q^2 + 1}\over 2} {\sf d} x_+ ~a_0 + {\bf i} 
                    {{q^2 - 1}\over 2} {\sf d}x_+ ~a_1, $$
$$a_1 ~{\sf d} x_+ = {{q^2 + 1}\over 2} {\sf d} x_+ ~a_1 - {\bf i} 
                    {{q^2 - 1}\over 2} {\sf d} x_+ ~a_0, $$
$$a_0 ~{\sf d} x_- = {{q^2 + 1}\over 2} {\sf d} x_- ~a_0 - {\bf i} 
                    {{q^2 - 1}\over 2} {\sf d} x_- ~a_1 + 
                    {({q - q^{-1})^2}\over 2} {\sf d} x_+ ~x_- - 
                    (q - q^{-1})({\sf d} a_2 ~a_2 + {\sf d} a_3 ~a_3),$$
$$a_1 ~{\sf d} x_- = {{q^2 + 1}\over 2} {\sf d} x_- ~a_1 + {\bf i} 
                    {{q^2 - 1}\over 2} {\sf d} x_- ~a_1 - {\bf i} 
                    {({q - q^{-1})^2}\over 2} {\sf d} x_+ ~x_- + {\bf i} 
                    (q - q^{-1})({\sf d} a_2 a_2 + {\sf d} a_3 ~a_3),$$
$$a_0 ~{\sf d} a_2 = q {\sf d} a_2 ~a_0 + {{q^2 - 1}\over 2} 
                      ({\sf d} a_0 + {\bf i} {\sf d} a_1) a_2, $$
$$a_0 ~{\sf d} a_3 = q {\sf d} a_3 ~a_0 + {{q^2 - 1}\over 2} 
                      ({\sf d} a_0 + {\bf i} {\sf d} a_1) a_3, $$
$$a_1 ~{\sf d} a_2 = q {\sf d} a_2 ~a_1 - {\bf i} {{q^2 - 1}\over 2} 
                      ({\sf d} a_0 + {\bf i} {\sf d} a_1) a_2, $$
$$a_1 ~{\sf d} a_3 = q {\sf d} a_3 ~a_1 - {\bf i} {{q^2 - 1}\over 2} 
                      ({\sf d} a_0 + {\bf i} {\sf d} a_1) a_3, $$
$$a_2 ~{\sf d} a_0 = q {\sf d} a_0 ~a_2 + {{q^2 - 1}\over 2} 
                      {\sf d} a_2 (a_0 - {\bf i} a_1), \eqno(24)$$
$$a_2 ~{\sf d} a_1 = q {\sf d} a_1 ~a_0 - {\bf i} {{q^2 - 1}\over 2} 
                      {\sf d} a_2 (a_0 - {\bf i} a_1), $$
$$a_2 ~{\sf d} a_2 = {{q^2 + 1}\over 2} {\sf d} a_2 ~a_2 - 
                    {{q^2 - 1}\over 2} {\sf d} a_3 ~a_3 - 
                    {{q - q^{-1}}\over 2} {\sf d} x_+ ~x_-, $$
$$a_2 ~{\sf d} a_3 = {{q^2 + 1}\over 2} {\sf d} a_3 ~a_2 + 
                    {{q^2 - 1}\over 2} {\sf d} a_2 ~a_3, $$
$$a_3 ~{\sf d} a_0 = q {\sf d} a_0 ~a_3 + {{q^2 - 1}\over 2} 
                      {\sf d} a_3 (a_0 - {\bf i} a_1), $$
$$a_3 ~{\sf d} a_1 = q {\sf d} a_1 ~a_3 - {\bf i} {{q^2 - 1}\over 2} 
                      {\sf d} a_3 (a_0 - {\bf i} a_1), $$
$$a_3 ~{\sf d} a_2 = {{q^2 + 1}\over 2} {\sf d} a_2 ~a_3 + {{q^2 - 1}\over 2} 
                      {\sf d} a_3 ~a_2, $$
$$a_3 ~{\sf d} a_3 = {{q^2 + 1}\over 2} {\sf d} a_3 ~a_3 - 
                    {{q^2 - 1}\over 2} {\sf d} a_2 ~a_2 - 
                    {{q - q^{-1}}\over 2} {\sf d} x_+ ~x_-, $$
where
$${\sf d} x_\pm = {\sf d} a_0 \pm {\bf i} {\sf d} a_1.$$

Applying the exterior differential {\sf d} on the relations (24) and using 
the nilpotency of {\sf d} one obtains 
$${\sf d} a_0 ~{\sf d} a_1 = - {\sf d} a_1 {\sf d} a_0, $$
$$({\sf d} a_0)^2 = 0 = ({\sf d} a_1)^2,$$
$${\sf d} a_0 ~{\sf d} a_2 = - {{q + q^{-1}}\over 2} {\sf d}a_2 ~{\sf d} a_0 
  + {\bf i} {{q - q^{-1}}\over 2} {\sf d} a_2 ~{\sf d} a_1,$$
$${\sf d} a_0 ~{\sf d} a_3 = - {{q + q^{-1}}\over 2} {\sf d} a_3 ~{\sf d} a_0 + 
          {\bf i} {{q - q^{-1}}\over 2} {\sf d} a_3 ~{\sf d} a_0,$$
$${\sf d} a_1 ~{\sf d} a_2 = - {{q + q^{-1}}\over 2} {\sf d} a_2 ~{\sf d} a_1 - 
          {\bf i} {{q - q^{-1}}\over 2} {\sf d} a_2 ~{\sf d} a_0,\eqno(25)$$
$${\sf d} a_1 ~{\sf d} a_3 = - {{q + q^{-1}}\over 2} {\sf d} a_3 ~{\sf d} a_1 - 
          {\bf i} {{q - q^{-1}}\over 2} {\sf d} a_3 ~{\sf d} a_1,$$
$${\sf d} a_2 ~{\sf d} a_3 = - {\sf d} a_3 ~{\sf d} a_2, $$
$$({\sf d} a_2)^2 = {\bf i} (q - q^{-1}) {\sf d} a_1 ~{\sf d} a_0 = 
                     ({\sf d} a_3)^2.$$

There is an interesting case which gives rise to the second kind of 
quaternionic variables, Grassmann quaternion is defined by 
$$\psi = \psi_0 e_0 + \psi_1 e_1 + \psi_2 e_2 + \psi_3 e_3, $$
where components $\psi_k$ $(k=0,1,2,3)$ are Grassmann variables. Essentially 
the relations (25) are the relations between the components $\psi_k$, 
in $q$-deformation. More details will be given in Appendix. 

To complete the differential calculus, we need the Cartan-Maurer one-forms. 
In analogy with the one-forms on a Lie group in classical differential 
geometry, one can construct the one-form $\Omega$ where 
$$\Omega = {\sf d} h~ h^\star = 
  w_0 e_0 + w_1 e_1 + w_2 e_2 + w_3 e_3. $$
So we can write the one-forms as follows 
$$w_0 = \alpha_0 a_0 + \alpha_1 a_1 + 
  {{q + q^{-1}}\over 2} (\alpha_2 a_2 + \alpha_3 a_3) + 
  {\bf i} {{q - q^{-1}}\over 2} (\alpha_3 a_2 - \alpha_2 a_3), $$
$$w_1 = - \alpha_0 a_1 + \alpha_1 a_0 - 
  {{q + q^{-1}}\over 2} (\alpha_2 a_3 - \alpha_3 a_2) - 
  {\bf i} {{q - q^{-1}}\over 2} (\alpha_2 a_2 + \alpha_3 a_3), $$
$$w_2 = \alpha_2 a_0 - \alpha_3 a_1 + 
  {{q + q^{-1}}\over 2} (\alpha_1 a_3 - \alpha_0 a_2) + 
  {\bf i} {{q - q^{-1}}\over 2} (\alpha_0 a_3 + \alpha_1 a_2), \eqno(26)$$
$$w_3 = \alpha_2 a_1 + \alpha_3 a_0 - 
  {{q + q^{-1}}\over 2} (\alpha_0 a_3 + \alpha_1 a_2) + 
  {\bf i} {{q - q^{-1}}\over 2} (\alpha_1 a_3 - \alpha_0 a_2), $$
where $\alpha_0 = {\sf d} a_0$, etc. 

We wish to find the commutation relations of the generators of ${\cal A}$ 
with those of the components of $\Omega$ which may be computed directly, 
as follows: 
$$a_0 w_+ = {{q^2 + 1}\over 2} w_+ a_0 + {\bf i} {{q^2 - 1}\over 2} w_+ a_1, 
  \qquad a_2 w_+ = {{q^2 + 1}\over 2} w_+ a_2 + {\bf i} {{q^2 - 1}\over 2} w_+ a_3,$$
$$a_1 w_+ = {{q^2 + 1}\over 2} w_+ a_1 - {\bf i} {{q^2 - 1}\over 2} w_+ a_0, 
  \qquad a_3 w_+ = {{q^2 + 1}\over 2} w_+ a_3 - {\bf i} {{q^2 - 1}\over 2} w_+ a_2,$$
$$a_0 w_2 = q w_2 a_0 + {{q - q^{-1}}\over 2} w_+ a_2, \qquad 
  a_2 w_2 = q w_2 a_2 - {{q - q^{-1}}\over 2} w_+ a_0,$$
$$a_1 w_2 = q w_2 a_1 - {{q - q^{-1}}\over 2} w_+ a_3, \qquad 
  a_3 w_2 = q w_2 a_3 + {{q - q^{-1}}\over 2} w_+ a_1, \eqno(27)$$
$$a_0 w_3 = q w_3 a_0 + {{q - q^{-1}}\over 2} w_+ a_3, \qquad 
  a_2 w_3 = q w_3 a_2 - {{q - q^{-1}}\over 2} w_+ a_1,$$
$$a_1 w_3 = q w_3 a_1 + {{q - q^{-1}}\over 2} w_+ a_2, \qquad 
  a_3 w_3 = q w_3 a_3 - {{q - q^{-1}}\over 2} w_+ a_0,$$
$$a_0 w_- = {{q^2 + 1}\over 2} w_- a_0 - {\bf i} {{q^2 - 1}\over 2} w_- a_1 + 
  {(q - q^{-1})^2 \over 2} w_+ (a_0 + {\bf i} a_1) + 
  (1 - q^2) (w_2 a_2 + w_3 a_3),$$
$$a_1 w_- = {{q^2 + 1}\over 2} w_- a_1 + {\bf i} {{q^2 - 1}\over 2} w_- a_0 - 
  {\bf i} {(q - q^{-1})^2 \over 2} w_+ (a_0 + {\bf i} a_1) - 
  (1 - q^2) (w_2 a_3 - w_3 a_2),$$
$$a_2 w_- = {{q^2 + 1}\over 2} w_- a_2 - {\bf i} {{q^2 - 1}\over 2} w_- a_3 + 
  {(q - q^{-1})^2 \over 2} w_+ (a_2 + {\bf i} a_3) + 
  (q^2 - 1) (w_2 a_0 + w_3 a_1),$$
$$a_3 w_- = {{q^2 + 1}\over 2} w_- a_3 + {\bf i} {{q^2 - 1}\over 2} w_- a_2 - 
  {\bf i} {(q - q^{-1})^2 \over 2} w_+ (a_2 + {\bf i} a_3) + 
  (1 - q^2) (w_2 a_1 - w_3 a_0),$$
where 
$$w_\pm = w_0 \pm {\bf i} w_1.$$
We now obtain the commutation relations of the Cartan-Maurer forms 
$$w_0^2 = {\bf i} {(q - q^{-1})^2\over 2} w_3 w_2, \qquad 
  w_1^2 = {\bf i} {{q^2 - q^{-2}}\over 2} w_3 w_2, $$
$$w_0 w_1 = - w_1 w_0 + (q^{-2} - 1) w_2 w_3,$$
$$w_0 w_2 = - w_2 w_0 + {{q^{-2} - q^2}\over 2} w_3 w_1 + {\bf i} 
  {(q - q^{-1})^2\over 2} w_2 w_1, $$
$$w_0 w_3 = - w_3 w_0 + {{q^2 - q^{-2}}\over 2} w_2 w_1 + 
  {\bf i} {(q - q^{-1})^2\over 2} w_3 w_1, $$
$$w_1 w_2 = - {{q^2 + q^{-2}}\over 2} w_2 w_1 + {\bf i} 
  {{q^{-2} - q^2}\over 2} w_3 w_1, \eqno(28)$$
$$w_1 w_3 = - {{q^2 + q^{-2}}\over 2} w_3 w_1 + 
  {\bf i} {{q^2 - q^{-2}}\over 2} w_2 w_1, $$
$$w_2 w_3 = - w_3 w_2, \qquad w_2^2 = 0 = w_3^2. $$

The conjugation of one-forms is defined  by 
$$w_0^\star = q^{-2} w_0 + {\bf i} (q^{-2} - 1) w_1, \qquad w_1^\star = w_1,$$
$$w_2^\star =  w_2, \qquad w_3^\star = w_3. \eqno(29)$$

Note that, we have 
$$\Omega + \overline{\Omega} = (1 - q^{-2}) (w_0 + {\bf i} w_1) e_0. \eqno(30)$$


To obtain the quantum algebra of the algebra generators we first write the 
Cartan- Maurer forms as 
$${\sf d} a_0 = w_0 a_0 - (w_1 a_1 + w_2 a_2 + w_3 a_3),$$
$${\sf d} a_1 = w_0 a_1 + w_1 a_0 + w_2 a_3 - w_3 a_2,$$
$${\sf d} a_2 = w_0 a_2 - w_1 a_3 + w_2 a_0 + w_3 a_1, \eqno(31)$$
$${\sf d} a_3 = w_0 a_3 + w_1 a_2 - w_2 a_1 + w_3 a_0.$$
Note that 
$${\sf d}^\star = q^2 {\sf d}. \eqno(32)$$
Indeed, for example, 
\begin{eqnarray*}
  ({\sf d} a_0)^\star 
& = & a_0^\star w_0^\star  - a_1^\star w_1^\star - a_2^\star w_2^\star - 
      a_3^\star w_3^\star \\
& = & q^{-2} a_0 (w_0 + {\bf i} w_1) - {\bf i} (a_0 - {\bf i} a_1) w_1 - 
      {1\over 2} a_2 [(q + q^{-1}) w_2 - {\bf i} (q - q^{-1}) w_3] \\ 
&  &  + {{\bf i}\over 2} a_3 [(q - q^{-1}) w_2 + {\bf i} (q + q^{-1}) w_3] \\
& = & q^2 (w_0 a_0 - w_1 a_1 - w_2 a_2 - w_3 a_3)
\end{eqnarray*}
so that 
$${\sf d}^\star a_0^\star = q^2 (w_0 a_0 - w_1 a_1 - w_2 a_2 - w_3 a_3) 
  = q^2 {\sf d} a_0 = q^2 {\sf d} a^\star_0 $$ 
implies that 
$${\sf d}^\star - q^2 {\sf d} = 0.$$

Using the nilpotency of the differential {\sf d}, we can write the two-forms 
as 
$${\sf d} w_0 = {\bf i} (q^{-2} - 1) w_2 w_3, \qquad 
  {\sf d} w_1 = (q^{-2} + 1) w_2 w_3,$$
$${\sf d} w_2 = {\bf i} (q^{-2} - 1) w_2 w_1 + (q^{-2} + 1) w_3 w_1, \eqno(33)$$
$${\sf d} w_3 = {\bf i} (q^{-2} - 1) w_3 w_1 - (q^{-2} + 1) w_2 w_1.$$

Using the Cartan-Maurer equations with together (27) we find the following 
commutation relations for the quantum algebra: 
$$\nabla_0 \nabla_1 = \nabla_1 \nabla_0, \qquad 
  \nabla_0 \nabla_2 = \nabla_2 \nabla_0, \qquad 
  \nabla_0 \nabla_3 = \nabla_3 \nabla_0, $$
\begin{eqnarray*}
\nabla_1 \nabla_2 
& = & {{q^2 + q^{-2}}\over 2} \nabla_2 \nabla_1 - 
    (q^{-2} + 1) \nabla_3 + {\bf i} (q^{-2} - 1) \nabla_2 - {\bf i} 
    {(q - q^{-1})^2\over 2} \nabla_0 \nabla_1 \\
&  & - {{q^2 - q^{-2}}\over 2} \nabla_3 (\nabla_0 + {\bf i} \nabla_1),  
\end{eqnarray*}
\begin{eqnarray*}
\nabla_1 \nabla_3 
& = & {{q^2 + q^{-2}}\over 2} \nabla_3 \nabla_1 + (q^{-2} + 1) \nabla_2 - 
      {\bf i} (q^{-2} - 1) \nabla_3 - 
      {\bf i} {(q - q^{-1})^2\over 2} \nabla_0 \nabla_1 \\
&  & + {{q^2 - q^{-2}}\over 2} \nabla_2 (\nabla_0 + {\bf i} \nabla_1), 
\end{eqnarray*}
\begin{eqnarray*}
\nabla_3 \nabla_2 
& = & \nabla_2 \nabla_3 + (q^{-2} + 1) \nabla_1 - {\bf i} (q^{-2} - 1) \nabla_0 
      + {\bf i} {{q^2 - q^{-2}}\over 2} \nabla_0^2 \\ 
&  & + {(q - q^{-1})^2\over 2} \nabla_1^2 + (1 - q^{-2}) \nabla_0 \nabla_1. 
\end{eqnarray*}

{\bf Acknowledgement}

This work was supported in part by {\bf T. B. T. A. K.} the 
Turkish Scientific and Technical Research Council. 

{\bf Appendix}

{\bf $q$-Deformation of Grassmann quaternionic algebra}

Quaternionic Lie groups can be defined in terms of matrices with quaternionic 
elements. To make contact with the usual formulation of Lie groups in terms 
of matrices with complex matrix elements, it is useful to represent each 
quaternion in terms of 2x2 complex matrices. Therefore, the quantum 
deformation of the Grassmann quaternionic algebra may be introduced using the 
idea of the quantum matrix theory [2,8]. 

Quantum Grassmann quaternionic algebra $\hat{H}_q$ is defined as a pair 
$(\hat{\cal A}, {\cal M})$, where $\hat{\cal A}$ is an algebra and 
has the following properties: 
$$\psi_0^2 = 0 = \psi_1^2, \qquad \psi_0 \psi_1 + \psi_1 \psi_0 = 0, $$
$$\psi_0 \psi_2 = - {{q + q^{-1}}\over 2} \psi_2 \psi_0 + 
                    {\bf i} {{q - q^{-1}}\over 2} \psi_2 \psi_1,$$
$$\psi_0 \psi_3 = - {{q + q^{-1}}\over 2} \psi_3 \psi_0 + 
                    {\bf i} {{q - q^{-1}}\over 2} \psi_3 \psi_1,$$
$$\psi_1 \psi_2 = - {{q + q^{-1}}\over 2} \psi_2 \psi_1 - 
                    {\bf i} {{q - q^{-1}}\over 2} \psi_2 \psi_0,$$
$$\psi_1 \psi_3 = - {{q + q^{-1}}\over 2} \psi_3 \psi_1 - 
                    {\bf i} {{q - q^{-1}}\over 2} \psi_3 \psi_0,$$
$$\psi_2^2 = {\bf i} {{q - q^{-1}}\over 2} \psi_1 \psi_0 = \psi_3^2, \qquad
  \psi_2 \psi_3 + \psi_3 \psi_2 = 0.$$

\vfill\eject

\end{document}